\documentclass[11pt]{amsart}
\usepackage{amsfonts}
\usepackage{amsmath}
\usepackage{fancyhdr,amssymb}
\usepackage{amsthm}
\usepackage{arydshln}
\usepackage{comment}
\usepackage{ color, ulem}
\input xy
\xyoption {all}

\input xy
\xyoption {all}

\def \Q {{\mathsf {Quot}}}
\def \P {{\mathsf P}}
\def \M {{\mathsf M}}
\def \e {{\mathcal E}}
\def \f {{\mathcal F}}

\def \V {{\text V}}
\DeclareMathOperator{\ch}{ch}

\newtheorem{theorem}{Theorem}

\theoremstyle{definition}
 
\theoremstyle{definition}

\theoremstyle{definition}

\newtheoremstyle{TheoremNum}
        {7pt}{7pt}              
        {\itshape}                      
        {}                              
        {\bfseries}                     
        {.}                             
        { }                             
        {\thmname{#1}\thmnote{ \bfseries #3}}
    \theoremstyle{TheoremNum}
    \newtheorem{thmn}{Theorem}

\begin{document}

\baselineskip=17pt

\title[The Segre-Verlinde correspondence]
{The Segre-Verlinde correspondence for the moduli space of stable bundles on a curve}
\author{Alina Marian}
\address{The Abdus Salam International Centre for Theoretical Physics, Strada Costiera 11, 34151 Trieste, Italy}
\address{Department of Mathematics, Northeastern University, 360 Huntington Avenue, Boston, MA 02115, USA }
\email{amarian@ictp.it}
\begin{abstract}
We show that the classic Verlinde numbers on the moduli space $\M(r,d)$ of rank $r$ and degree $d$ semistable vector bundles over a smooth projective curve can also be regarded as Segre numbers of natural universal complexes over $\M(r,d).$ This leads to interesting identities among universal integrals on $\M(r,d).$
\end{abstract}
\maketitle

\vskip.3in

\section{The Segre-Verlinde identity}

\vskip.2in

Let $C$ be a smooth projective curve of genus $g \geq 2.$ The classic Mumford moduli space $\M (r,d)$ of rank $r$ degree $d$ slope-semistable vector bundles on $C$ is a projective variety of dimension $r^2 (g-1) +1.$ Its Picard group is generated by determinant line bundles defined relative to vector bundles on $C$ of slope $\mu_0 = g -1 - \frac{d}{r}$ and minimal possible rank (cf. \cite{dn}). More precisely, to such a reference bundle $F_0 \to C,$ one associates a line bundle $$\Theta_{r, F_0} \in \text{Pic} \, \M (r, d),$$ constructed by GIT descent; for any family $\mathcal V \to T \times C$ of semistable rank $r$ and degree $d$ vector bundles over a smooth base $T$, with associated morphism $f_T: T \to \M (r, d)$, we have $$f_T^\star \Theta_{r, F_0} = \det \mathsf R\pi_\star (\mathcal V \otimes \rho^\star F_0 )^{-1}.$$ Here $\pi, \rho$ denote the projections from the product $T \times C$ to the two factors. For further use we also set 
$$\det: \M (r, d) \to \text{Jac}^d (C), \, \, E \to \det E,$$ the natural map sending bundles to their determinants.

To state the observation of this paper in simplest form, let us first restrict to the case when $r$ and $d$ are coprime. Then the moduli space $\M (r, d)$ is smooth, parametrizes only stable bundles, and admits a universal vector bundle $$\mathcal {V} \to \M (r,d) \times C.$$   We show 
\begin{theorem}
Assume $\text{gcd}\, (r, d) = 1,$ let $\ell \geq 1$, and let $\alpha \in K (C)$ be a K-class with
$$\text{rank} \, \, \alpha = (\ell +1) r , \, \, \, \text{degree} \,\,  \alpha = - (\ell +1) d + \ell r (g-1).$$ 
Set $$\alpha_{\M} = \mathsf R\pi_\star (\mathcal V \otimes \rho^\star \alpha),$$ a $K$-class on $\M (r, d).$
Then 
\begin{equation}
\label{segverlinde}
\chi (\M(r, d), \, \Theta_r^\ell \otimes {\det}^\star \Theta_1) = \int_{\M (r, d)} s (\alpha_{\M}).
\end{equation}
\label{coprime}
\end{theorem}

The two sides of \eqref{segverlinde} have different flavor. The left side is the classic Verlinde number giving the holomorphic Euler characteristic of the determinant line bundle $\Theta_{r, F_0}$ at level $\ell$, twisted by the theta line bundle $\Theta_1$ from the Jacobian of $C$. The reference minimal-rank bundle $F_0 \to C$ as well as a suitable reference line bundle for $\Theta_1$ are suppressed from the notation for simplicity, since the Euler characteristic is independent of their choice. 

The right side is the Segre class of $\alpha_\M \in K (\M (r, d))$. Note that 
\begin{eqnarray*}
\text{rank}\, \,  \alpha_{\M} &= & \chi (C, E \cdot \alpha) \, \, \text{for} \, \, E \in \M (r, d) \\
& = & r d (\ell +1) + r [  - (\ell +1) d + \ell r (g-1)] - r^2 (\ell +1) (g-1) \\ &=& -r^2 (g-1).
\end{eqnarray*}
The push-forward complex $\alpha_\M$ is not invariant under tensoring the universal bundle $\mathcal V \to \M (r, d) \times C$ by a line bundle on $\M (r,d),$ but its rank guarantees that its top Segre class is in fact invariant.
Indeed, for any variety $X$, line bundle $L \to X$, and $K$-class $W$ on $X$ of rank $n$, it is immediate that
$$c_{n+1} (W \otimes L ) = c_{n+1} (W).$$ We have $$ \int_{\M (r, d)} s (\alpha_{\M}) =  \int_{\M (r, d)} c_{r^2 (g-1) + 1} (- \alpha_{\M}),$$ and the latter class is invariant under twisting the complex $\alpha_\M$ by a line bundle. 

Furthermore, both sides are polynomials of degree $r^2 (g-1) + 1$ in the level $\ell$. It is easy to see, for both Segre and Verlinde numbers, that the coefficient of the highest order term in $\ell$ recovers the volume of the moduli space $\M (r,d)$. Matching the remaining coefficients yields however unexpected identities between tautological top intersections on $\M(r,d).$ 

\vskip.2in

All is not lost for equality \eqref{segverlinde} in the presence of strictly semistables, for an arbitrary pair $(r, d).$ We record the general slightly less direct statement as a variation on Theorem \ref{coprime}. 
\begin{thmn}[\ref{coprime}$^*$]
\label{arbitrary_rd}
For any rank-degree pair $(r, d)$, set $h = \text{gcd}\,  (r, d), \, \, \, r_0 = \frac{r}{h}, \, \, \, d_0 = \frac{d}{h}.$ Let $\ell \geq 1,$ and fix a point $p \in C$. Let $\alpha \in K (C)$ be the K-class with
$$\text{rank} \, \, \alpha = (\ell +h) r_0, \, \, \, \, c_1 (\alpha)  = (- (\ell +h) d_0 + \ell r_0 (g-1)) \,[p].$$ There exists a zero-dimensional class $s(\alpha_\M ) \in A_0 (\M (r, d))$ so that

(i) For any family $\mathcal V \to T \times C$  of {\it stable}  rank $r$ and degree $d$ vector bundles on $C$ with a flat classifying map $f_T: T \to \M(r,d),$ we have 
$$f_T^{\star} s (\alpha_\M) = s_{r^2 (g-1) + 1} \left (\mathsf R\pi_\star ( \mathcal V \cdot \rho^{\star} \alpha ) \right ).$$

(ii) $\chi (\M(r, d), \, \Theta_r^\ell \otimes {\det}^\star \Theta_1) = \int_{\M (r, d)} s (\alpha_{\M}).$

\vskip.1in
\noindent
In the distinguished case $(r, 0)$ we simply have $$\text{rank} \, \,\alpha = \ell +r, \, \, \, \, \text{degree} \, \, \alpha =  \ell  (g-1).$$

\end{thmn} 

\vskip.2in

While the geometry and topology of $\M (r,d)$ have been extensively explored across decades, it is not easy to probe equality \eqref{segverlinde} directly for an arbitrary pair $(r, d)$. Nevertheless, the structure and symmetries of important invariants of $\M (r,d)$
have long been understood in connection with the Grothendieck Quot scheme  $\Q (\mathbb C^n, r, d)$ parametrizing rank $r$ degree $-d$ subsheaves $$E \hookrightarrow \mathbb C^n \otimes {\mathcal O}_C \, \, \, \text{on} \, \, \, C.$$The beautiful representation of the Verlinde number as a distinguished intersection on the Quot scheme, due to Witten \cite{witten}, is a first step in the argument toward Theorem \ref{coprime}. The second is the correction-free wall-crossing formula of \cite{alina} allowing one to move the Quot scheme intersection back to the moduli space $\M (r, d)$, where the Segre class emerges. 
The two ingredients are substantial but available, so Theorem \ref{coprime} is established quickly. The context is described over the next two sections. The wall-crossing is stated as Theorem \ref{wallcrossing} below. A few explicit tautological identities resulting from Theorem \ref{coprime} are written in Section \ref{examples}. The last section reflects on the resulting triangle of invariants and places it in the broader setting of recent developments for Segre-Verlinde dualities. 

\vskip.4in

\section{The Quot scheme and Witten's intersection} 

\vskip.2in

We consider the Grothendieck Quot scheme  $\Q (\V, r, d)$ parametrizing rank $r$ degree $-d$ subsheaves $$E \hookrightarrow \V \, \, \, \text{on} \, \, \, C,$$ for a fixed vector bundle $\V \to C$ of rank $n$ and degree $\delta$. 
The Quot scheme carries a universal family
$$0 \to \e \to \rho^\star \V \to \f \to 0 \, \, \, \text{on} \, \, \, \Q (\V, r, d) \times C,$$ and a virtual class $$[\Q (\V, r, d)]^{vir} \in A_\star (\Q (\V, r, d))$$ of dimension 
\begin{equation}
 \text{vdim}\, \, \Q (\V, r, d) = n \,d  + r \delta - r (n-r) (g-1),
 \label{expecteddim}
 \end{equation}
 constructed in \cite{mo1}. 
The virtual intersection theory of the universal Chern classes $$c_i (\e^\vee) \in A^i (\Q (\V, r,d) \times C), \, \, \, 1\leq i \leq r,$$ has a rich structure. For brevity we refer to the various slant products of the universal $c_i (\e^\vee)$ with classes on $C$ as {\it Atiyah-Bott classes}, connecting back to the classic \cite{ab}. 

\vskip.2in

 In the basic case $\V = \mathbb C^n,$ $\Q (\mathbb C^n, r, d)$ is a compactification of the space of maps of degree $d$ from $C$ to the Grassmannian $G (r, n).$ For large degree $d$, the Quot scheme is moreover irreducible \cite{bdw}.
In this context, the top intersections of the classes $$a_i = c_i (\e_p^\vee) \in A^i (\Q (\mathbb C^n, r, d)),\, \, \text{for}\, \, p \in C,$$ were studied in \cite{bertram}, and were shown to have enumerative meaning: they calculate counts of maps from $C$ to the Grassmannian $G (r,n)$ sending fixed points of $C$ to special Schubert subvarieties of $G (r, n).$ 
Furthermore, these top intersections are captured by a beautiful closed formula, the Vafa-Intriligator formula proposed in \cite{i} and proven in \cite{st}, \cite{bertram}, \cite{r}, \cite{mo1}. 

\vskip.2in

Amidst all, the Atiyah-Bott class $$a_r = c_r (\mathcal E_p^\vee) \in A^r (\Q (\V, r, d)), \, \, \text{for} \, \, p \in C,$$ plays a distinguished role. An elementary modification $0 \to \V_0 \to \V \to \mathbb C_p \to 0$ at $p$ gives an inclusion
$$\iota: \Q (\V_0, r, d) \to \Q (\V, r, d),$$ under which we have
$$\iota_\star [\Q (\V_0, r, d)]^{vir} = a_r \cdot [\Q (\V, r, d)]^{vir} \in A_\star (\Q (\V, r, d)).$$
Furthermore, the top self-intersection of $a_r$ on the Quot scheme associated with $$\V = \mathbb C^{\ell +r}$$ calculates the Verlinde theory at level $\ell$ on the moduli space $\M (r, 0)$ of semistable bundles of rank $r$ and degree $0$ on $C$.
Choosing $d$ divisible by both $r$ and $\ell$ we have
\begin{equation}
\chi (\M (r, 0), \, \Theta^\ell_r \otimes \Theta_1)=  \int_{[\Q (\mathbb C^{\ell +r }, r, d)]^{vir}} a_r^{\ell \, t},
\label{witten0}
\end{equation}
where $t$ is the level-rank symmetric quantity
$$ t = \frac{d}{r} + \frac{d}{\ell} - (g-1).$$

Formula \eqref{witten0} is one of the most beautiful identities in enumerative geometry.  On the Grassmannian $G (r, \ell +r)$ with tautological subbundle $S$, the top intersection $c_r (S^\vee)^\ell$ represents the class of a point, 
$$c_r (S^\vee)^\ell = [\text{point}] \in A_0 (G (r, \ell + r)).$$
Thus for sufficiently large $d$ divisible by both $r$ and $\ell$, the identity \eqref{witten0} realizes the Verlinde number at level $\ell$ simply as the count of maps from $C$ to the Grassmannian $G (r, \ell + r)$ sending $t$ fixed points $p_1, \ldots, p_t$ on C to fixed points $q_1, \ldots, q_t$ in $G (r, \ell + r).$ This number is moreover calculated by a special case of the Vafa-Intriligator formula. The equality \eqref{witten0} is due to Witten \cite{witten}, and was formulated on the Quot scheme using the Atiyah-Bott class $a_r$ in \cite{mo2}. The formula makes clear the symmetry of the Verlinde numbers under exchanging rank and level, and also enhances it.  Indeed, beyond the numerical level-rank symmetry,  the formula is consistent with a {\it geometric} strange duality isomorphism involving the spaces of sections of the line bundles with dimension given by the Verlinde numbers \eqref{witten0}, as explained in \cite{mo2}.   

\vskip.1in

For an arbitrary pair $(r,d),$ we set as before $$h = \text{gcd}\, (r, d), \, \, \, r_0 = \frac{r}{h}, \, \, \, d_0 = \frac{d}{h}. $$ Since $\M (r, d)$ only depends on $d$ modulo $r$, we may assume $d$ is suitably large. Generalizing \eqref{witten0}, we can write (cf. \cite{mo2}, Section 5)

\begin{equation}
\label{witten1} 
 \chi (\M (r, d), \, \Theta^\ell_r \otimes \Theta_1) = \int_{[\Q (\mathbb C^{(\ell + h) r_0}, r, d)]^{vir}} a_r^{(\ell + h) d_0 - \ell r_0 (g-1)}.
\end{equation}
This formula is the first of the two ingredients needed to establish Theorem \ref{coprime}. 

\vskip.4in

\section{The moduli space of stable pairs and the wall-crossing formula} 

\vskip.2in

Turning to the second ingredient, by taking duals of exact sequences $$0 \to E \to {\mathcal O}^n \to F \to 0 \, \, \, \text{on} \, \, \, C,$$ we view the Quot scheme as the moduli space of pairs $$\{\tau: {{\mathcal O}^n}^\vee \to E^\vee, \,\,  \text{with} \, \, E^\vee \, \,  {\text{pure and dim coker}} \,\, \tau = 0\}.$$ Moduli spaces of stable pairs on $C$ were constructed in \cite{bradlow}, \cite{bd}, \cite{bdw}, depending on a stability parameter which sweeps an admissible range. They can also be viewed as a special case, in a curve setting, of the general algebro-geometric construction carried out in \cite{lin}. 

The Quot scheme $\Q (\mathbb C^n, r, d)$ is the large parameter limit of this construction. At the other end of the stability range we find the moduli space $$\P (\mathbb C^n, r, d) = \{ \tau: \mathcal O^n \to E, \,\,  \text{with} \,\,  E \, \, \text{slope semistable and} $$ $$\mu (F ) < \mu (E)\, \text{for any proper} \, F\subset E  \, \text{with Im} \, \tau \subset F\}.$$ For $d > 2r (g-1),$ this moduli space is a smooth projective variety endowed with a universal morphism
$$\tau: \mathcal O^n \to \mathcal U \, \, \, \text{on} \, \, \, \P (\mathbb C^n, r, d) \times C,$$ 
and a surjective map forgetting the sections,  $$\Psi: \P (\mathbb C^n, r, d) \to \M (r, d).$$ When $\text{gcd}\, (r,d) = 1$ one notices immediately that
\begin{equation}
 \label{stablepairs}
 \P (\mathbb C^n, r, d) = \mathbb P (\pi_\star \mathcal V^{\oplus n}),
\end{equation}
 a projective bundle over $M (r,d).$ 
The universal bundles match, 
\begin{equation}
\mathcal U = \mathcal O (1) \otimes (\Psi\times 1)^\star \mathcal V\, \, \, \text{on} \, \, \, \P (\mathbb C^n, r, d) \times C,
\end{equation}
where $\mathcal O (1) \to \P (\mathbb C^n, r, d)$ is the associated hyperplane line bundle. 
 
 \vskip.1in
 
 We set up the wall-crossing in the regime when the degree $d$ is large (relative to $r, n, g$), so the Quot scheme $\Q (\mathbb C^n, r, d)$ is irreducible. In this case, the two spaces $\Q (\mathbb C^n, r, d)$ and $\P (\mathbb C^n, r,d)$ are birational of dimension
$nd - r (n-r) (g-1),$
 agreeing on the open subscheme $$\Q^{ss} = \{\tau: \mathcal O^n \to E, \, \dim \text{coker} \, \tau = 0, \, \, E \, \text{semistable} \}.$$
 The two universal bundles $\mathcal E^\vee$ and $\mathcal U$ also agree here. The following equality holds. 
 \vskip.1in
 
 
\begin{theorem} [\cite{alina}] 
\label{wallcrossing}
Let $n > r$ and assume $d$ is sufficiently large relative to $r, n, g$.   Assume $r$ divides $nd$, letting $r m = nd - r (n-r) (g-1). $ 
We have
\begin{equation}
\label{nocorrections}
\int_{\P (\mathbb C^n, r, d)} c_r(\mathcal U_p)^m = \int_{\Q (\mathbb C^n, r, d)} c_r (\mathcal E_p^\vee) ^m.
\end{equation}
\end{theorem}

The remarkable fact expressed by the theorem is that no wall-crossing corrections are needed for the top intersection of the special Atiyah-Bott class. The typically arduous task of crossing the walls is thus avoided under the theorem. This phenomenon was studied more broadly in \cite{alina}, with the goal of calculating the full intersection theory of $\M (r,d)$ on the Quot scheme. Theorem 1 of \cite{alina} explains that any polynomial of a fixed degree $k$ in the tautological Atiyah-Bott classes paired by a complementary power of the distinguished cycle $a_r$ can be equally evaluated on the Quot scheme and on the moduli space of stable pairs, {\it provided the number $n$ of sections is suitably large relative to $k$}. As an example, the full intersection theory of the moduli space of stable rank 2 odd degree vector bundles on $C$ was evaluated on the Quot scheme in \cite{mo3}. Furthermore, the Verlinde numbers were calculated in generality in \cite{mo4} from the large $n$ asymptotics of the Vafa-Intriligator formula.

The argument in \cite{alina} relies on a careful examination of the transversality properties of the special cycle $a_r$ relative to the boundary loci of the two birational moduli spaces. If one is solely interested in evaluating the top power of the special cycle $a_r$, as is the case in Theorem \ref{wallcrossing}, the condition\begin{footnote} {The condition $n\geq 2r+1$ stated in Theorem 1 of \cite{alina} eases the intermediate statement of Lemma 1 therein but is not actually needed within the setup of Theorem 1. }\end{footnote} $n \geq 2r+1$ and the lower bound on the number $n$ of sections in \cite{alina} give way to the simple requirement $n > r$.

\vskip.2in

In the case when $\text{gcd} \,(r,d) = 1$ it is easy to further calculate the left-hand side of \eqref{nocorrections}, pushing it forward by the forgetful morphism $\Psi: \P (\mathbb C^n, r, d) \to \M (r, d).$ We set first $$\zeta = c_1 (\mathcal O (1)) \, \, \text{on}\, \, 
 \P (\mathbb C^n, r, d),$$
 and find:

\begin{eqnarray}
\int_{\P (\mathbb C^n, r, d)} c_r (\mathcal U_p)^m  &= & \int_{\P (\mathbb C^n, r, d)} c_{\text{top}} \left ( \mathcal U_p^{\, \oplus m} \right ) 
=  \int_{\P (\mathbb C^n, r, d)} c_{\text{top}} \left ( \mathcal O (1) \otimes ( \mathcal V_p^{\, \oplus m}) \right ) \nonumber \\
&= &\int_{\P (\mathbb C^n, r, d)} \frac{1}{1-\zeta} \cdot c \left ( \mathcal V_p^{\, \oplus m} \right ) \nonumber \\
&=& \int_{\M (r, d)} s (\pi_\star\mathcal V^{\,\oplus n}) c \left ( \mathcal V_p^{\, \oplus m} \right )\nonumber  \\
&= & \int_{\M (r, d)} s ({\mathsf R}\pi_\star ( \mathcal V\otimes (\mathbb C^n - \mathbb C_p^{\, \oplus m} ))) \nonumber \\
&=&  \int_{\M (r, d)} s (\alpha_\M), \, \, \text{for} \, \, \alpha \in K (C), \, \text{rank}\,  \alpha = n, \, \deg \alpha = -m. 
\label{segre}
\end{eqnarray}

\vskip.1in

For an arbitrary pair $(r, d)$, the pushforward class $$\Psi_\star \left ( c_r ( \mathcal U_p )^m \right ) \in A_0 (\M (r, d))$$ is well defined. For a family $\mathcal V \to T \times C$ of stable rank $r$ degree $d$ bundles on $C$ with a flat associated morphism $f_T: T \to \M (r, d)$, we have a corresponding map
$$f_{\mathbb P}: \mathbb P ( \pi_\star \mathcal V^{\, \oplus n} ) \to \P (\mathbb C^n, r, d).$$
Here $\Psi_T: \mathbb P ( \pi_\star \mathcal V^{\, \oplus n} ) \to T$ is the projective bundle parametrizing vector bundles in the family $T$ together with $n$ sections. Then  
\begin{equation} 
f_T^\star \Psi_\star \left ( c_r ( \mathcal U_p )^m \right ) = {\Psi_T}_\star f_{\mathbb P}^\star \left ( c_r ( \mathcal U_p )^m \right ) =  {\Psi_T}_\star  \left ( c_r ( \mathcal O (1) \otimes \mathcal V_p )^m \right ) = s_{r^2 (g-1) + 1} \left ( \mathsf R\pi_\star (\mathcal V \cdot \rho^\star\alpha) \right ),
\label{basechange}
\end{equation}
by the same calculation \eqref{segre}. 

\vskip.1in

Using the Quot scheme expression \eqref{witten1} of the general Verlinde number,  the wallcrossing equality \eqref{nocorrections} of Theorem \ref{wallcrossing}, and the simple calculations \eqref{segre}, \eqref{basechange}, Theorems \ref{coprime} and {\ref{coprime}}$^*$ follow. We have: 
\begin{eqnarray*}
 \chi (\M (r, d), \, \Theta^\ell_r \otimes \Theta_1) &=& \int_{[\Q (\mathbb C^{(\ell + h) r_0}, r, d)]^{vir}} a_r^{\,(\ell + h) d_0 - \ell r_0 (g-1)} \\
& = &  \int_{\P (\mathbb C^{(\ell + h) r_0}, r, d)} c_r (\mathcal U_p)^{\, (\ell + h) d_0 - \ell r_0 (g-1)} \\
& = &  \int_{\M (r, d)} s (\alpha_\M),  
\end{eqnarray*}
for $\alpha \in K (C), \, \text{rank}\,  \alpha = (\ell + h) r_0, \, \deg \alpha = - (\ell +h) d_0 + \ell r_0 (g-1).$
\qed

\vskip.4in

\section{Examples}
\label{examples}

\vskip.2in

We write down a few identities of tautological top intersections, resulting from Theorem 1. We assume here that $gcd(r, \, d) = 1.$ It is easier to express the equalities on the simply-connected moduli space $\mathsf N (r, d)$ of stable bundles with fixed determinant, which has dimension $(r^2 -1) (g-1).$ We let the universal bundle be 
$$\mathcal W \to \mathsf N (r, d) \times C.$$ We K\"{u}nneth-decompose its Chern classes with respect to a symplectic basis $\{1, \, \delta_j, \, \omega \},$ $ 1\leq j \leq 2g,$ for the cohomology of the curve $C$: 
\begin{equation}
\label{universal}
c_i (\mathcal W) = \mathsf a_i \otimes 1 + \sum_{j=0}^{2g} \mathsf b_i^j \otimes \delta_j + \mathsf f_i \otimes \omega, \, \, 1\leq i \leq r.
\end{equation}
It is well known (cf \cite{ab}) that the classes $\mathsf a_i, \mathsf b_i^j, \mathsf f_i$ where $i\leq r$ generate the cohomology ring $H^\star (\mathsf N (r,d), \, \mathbb Q)$ multiplicatively. 
For convenience, we formally normalize the universal $\mathcal W$ (by twisting it with $\det \mathcal W_p^{-\frac{1}{r}}$) so that 
\begin{equation}
\label{normalize}
\mathsf a_1 = 0.
\end{equation}
We then have
$$c_1 (\mathcal W) = d\omega, \, \, \, c_1 (\Theta_r) = r \,\mathsf f_2,$$
with $\text{Pic} \, \mathsf N (r, d) = \mathbb Z \,  \Theta_r.$ 
Letting $\mathsf J$ denote the Jacobian of degree 0 line bundles on $C$, we consider the \'{e}tale morphism 
$$\tau: \mathsf N (r, d) \times \mathsf J \to \mathsf M (r, d), \, \, \, \, (E, L) \to E \otimes L.$$
Note that $\tau$ has degree $r^{2g}$. By pulling back under $\tau$, we calculate
\begin{eqnarray*}
\chi (\mathsf M (r, d), \, \Theta_r^\ell \otimes \Theta_1) &=& \frac{1}{r^{2g}} \, \chi (\mathsf N (r, d), \,  \Theta_r^{\ell} ) \cdot \chi (\mathsf J, \,    \Theta_1^{r^2 \ell + r^2}) \\ &=&  \frac{1}{r^{2g}} (r^2 \ell + r^2)^g \cdot  \chi (\mathsf N (r, d), \, \Theta_r^\ell) \\ &=& (\ell +1)^g  \cdot \chi (\mathsf N (r, d), \,  \Theta_r^\ell).
\end{eqnarray*}
In turn,
\begin{equation}
\label{verlindeN}
\chi (\mathsf N (r, d), \,  \Theta_r^\ell) = \int_{\mathsf N (r, d)} \exp (\ell r \, \mathsf f_2) \cdot \text{Todd} \, \mathsf N (r, d) = \int_{\mathsf N (r, d)} \exp ((\ell +1)r \, \mathsf f_2) \cdot \hat{\mathsf A} (\mathsf N (r, d)),
\end{equation}
where the $\hat{\mathsf A}$ genus is a polynomial in the classes $\mathsf a_i = c_i (\mathcal W_p), 1\leq i \leq r,$ explicitly given in terms of the Chern roots $\alpha_i$ of $\mathcal W_p$ as
$$\hat{\mathsf A} (\mathsf N (r, d)) = \prod_{1\leq i<j \leq r} \left ( \frac{\alpha_i - \alpha_j}{2 \sinh \frac{\alpha_i - \alpha_j}{2}} \right )^{2 (g-1)}.$$
 
 \vskip.2in
 
 \subsection{High-degree terms in arbitrary rank}
 
 It is clear that $ \chi (\mathsf N (r, d), \,  \Theta_r^\ell)$ is a polynomial of degree $(r^2-1)(g-1)$ in $\ell +1$. Furthermore, since the $\hat{\mathsf A}$ genus is an even function, all terms of degree equal to $(r^2-1)(g-1) - 1 \mod 2$ vanish. Setting $$t = (r^2-1) (g-1),$$ we can calculate explicitly the highest-order terms occurring in the expansion of the Verlinde number $\chi (\mathsf M (r, d), \, \Theta_r^\ell \otimes \Theta_1)$:
\begin{eqnarray}
\chi (\mathsf M (r, d), \,  \Theta_r^\ell \otimes \Theta_1) &=& (\ell+1)^{t +g}  \int_{\mathsf N (r, d)} \frac{(r \, \mathsf f_2)^t}{t!} \label{verexp}\\
&+& (\ell+1)^{t+g -2}\, \, \frac{r (g-1)}{6}\int_{\mathsf N (r, d)} \frac{(r\, \mathsf f_2)^{t-2}}{(t -2)!} \cdot \mathsf a_2\nonumber  \\ 
& + & (\ell+1)^{t+g -4}  \int_{\mathsf N (r, d)} \frac{(r\, \mathsf f_2)^{t-4}}{(t -4)!} \cdot \left [ - \frac{r (g-1)}{360} \mathsf a_4 \, + \right. \nonumber \\
& + & \left. \left ( \frac{r^2 (g-1)^2}{72} + \frac{r (g-1)}{720}  + \frac{g-1}{120}   \right ) \mathsf a_2^2 \right ] \nonumber \\
& + & {\mathcal Order} (\ell +1)^{t+g -6}. \nonumber
\end{eqnarray}

All terms in the above Riemann-Roch polynomial are top intersections involving solely the classes $\mathsf f_2, \mathsf a_2, \, \ldots, \mathsf a_r.$ The situation is different for the evaluation of the Segre class of $\alpha_\mathsf M.$ We have
$$\int_{\mathsf M (r, d)} s (\alpha_{\mathsf M}) = \frac{1}{r^{2g}} \int_{\mathsf N (r, d) \times \mathsf J} s (\tau^\star \alpha_{\mathsf M}),$$
where $$\tau^\star \alpha_{\mathsf M} = \mathsf R\pi_\star ( \mathcal W \otimes \mathcal P \otimes \rho^\star \alpha),$$
with $\pi: \mathsf N (r, d) \times \mathsf J \times C \longrightarrow \mathsf N (r, d) \times \mathsf J$ denoting the projection to the first two factors. 
Here $\mathcal P \to \mathsf J \times C$ is the Poincar\'{e} line bundle, with $$\ch \mathcal P = 1 + \sum_{i=1}^{2g} \gamma_i \otimes \delta_i - \theta \cdot \omega \in H^\star (\mathsf J \times C),$$ 
where $$\theta = c_1 (\Theta_1) = \sum_{i=1}^g \gamma_i \gamma_{i+g} \in H^2 (\mathsf J, \, \mathbb Z)$$ denotes the principal polarization on the Jacobian. 

Applying the Grothendieck-Riemann-Roch theorem, we find
 $$\text{ch} \, \tau^\star \alpha_{\mathsf M} = \pi_\star \left [ \text{ch} \, \mathcal W \cdot \text{ch} \, \mathcal P \cdot \left ( (\ell +1) r  - d (\ell+1) \omega - r (g-1) \omega \right ) \right ].$$
Setting $$\alpha_{\mathsf N} = \mathsf R\pi_\star ( \mathcal W  \otimes \rho^\star \alpha) \in K (\mathsf N (r, d)), $$
we see that
\begin{equation}
\label{segreMN}
\text{ch} \, \tau^\star \alpha_{\mathsf M} = \text{ch} \, \alpha_{\mathsf N} - (\ell +1) r \cdot \theta\cdot  \text{ch} \, \mathcal W_p + (\ell + 1) r \cdot \pi_\star \left ( \text{ch} \, \mathcal W \cdot \sum_{i=1}^{2g} \gamma_i \otimes \delta_i \right ).
\end{equation}
Passing standardly from Chern character to Segre class, we find the highest-order terms in $\ell+1$ of the Segre polynomial to be 
\begin{eqnarray}
\int_{\mathsf M (r, d)} s (\alpha_{\mathsf M}) & = & (\ell +1)^{t+g}  \int_{\mathsf N (r, d)} \frac{(r \, \mathsf f_2)^t}{t!}  \\
&+& (\ell+1)^{t+g -1} \int_{\mathsf N (r, d)} \frac{(r\, \mathsf f_2)^{t-2}}{(t -2)!} \cdot \left [ \left( d-\frac{dr}{2} \right ) \mathsf a_2 + \frac{r}{2} \mathsf f_3 \right ] \label{segreexp1}\\
& + & {\mathcal Order}  (\ell +1)^{t+g-2}, \nonumber
\end{eqnarray} 
where we write down the next two terms explicitly as follows. 
The term with order $(\ell + 1)^{t+ g-2}$ equals
\begin{equation}
 \label{segreexp2}
 \int_{\mathsf N (r, d)} \frac{(r\, \mathsf f_2)^{t-2}}{(t -2)!} \cdot \mathsf a_2 \left ( -\frac{(r^2-1) (g-1) -2}{3} +  r (g-1) - \frac{2g}{r} \right) 
 \end{equation}
$$+  \int_{\mathsf N (r, d)} \frac{(r\, \mathsf f_2)^{t-3}}{(t -3)!} \cdot \left [ \left (d -\frac{dr}{3}\right ) \mathsf a_3 + \frac{r}{3} \mathsf f_4 + \sum_{j=1}^g \mathsf b_2^j \mathsf b_2^{j+g}  \left ( \frac{r}{3} - 2 \right )       \right ] $$
$$ +   \int_{\mathsf N (r, d)} \frac{(r\, \mathsf f_2)^{t-4}}{(t -4)!} \cdot \left [ \frac{r}{2} \mathsf f_3 + d \left ( 1 - \frac{r}{2} \right ) \mathsf a_2  \right ]^2.
$$
\vskip.2in

\noindent The term with order $(\ell + 1)^{t+ g-3}$ equals

\vskip.1in

\begin{equation} 
\label{segreexp3}
 \int_{\mathsf N (r, d)} \frac{(r\, \mathsf f_2)^{t-3}}{(t -3)!} \cdot \left ( -\frac{3g}{r} + r (g-1) \right ) \mathsf a_3 \, + 
 \end{equation}
$$+
\int_{\mathsf N (r, d)} \frac{(r\, \mathsf f_2)^{t-4}}{(t -4)!} \cdot \left [ \frac{r}{4} \mathsf f_5 + \left ( r^2 (g-1) - g - \frac{r}{4} \right )\mathsf a_2 \mathsf f_3 -  \frac{r}{4} \mathsf a_3 \mathsf f_2 + d \left ( 1 - \frac{r}{4} \right ) \mathsf a_4 +\right. $$
$$\left. + \, d \left (\frac{1}{2}\left ( \frac{r}{2} - 1 \right) + r (2-r) (g-1) - (2-r) \frac{g}{r} \right ) \mathsf a_2^2 +  \left ( 2 -\frac{r}{4} \right ) \sum_{j=1}^g ( \mathsf b_2^j \mathsf b_3^{j+g} + \mathsf b_3^j \mathsf b_2^{j+g} ) \right ] $$
$$ + \int_{\mathsf N (r, d)} \frac{(r\, \mathsf f_2)^{t-5}}{(t -5)!} \cdot \left [ \frac{r}{2} \mathsf f_3 + d \left (1 - \frac{r}{2} \right ) \mathsf a_2 \right ] \cdot \left [ d \left ( 1 - \frac{r}{3} \right ) \mathsf a_3   + \frac{r}{3} \mathsf f_4  - \frac{r}{3} \mathsf a_2 \mathsf f_2 + \frac{r}{3}  \sum_{j=1}^g  \mathsf b_2^j \mathsf b_2^{j+g} \right ] $$
$$ +   \int_{\mathsf N (r, d)} \frac{(r\, \mathsf f_2)^{t-6}}{(t -6)!} \cdot \left [ \frac{r}{2} \mathsf f_3 + d \left ( 1 - \frac{r}{2} \right ) \mathsf a_2\right ]^3.$$

\vskip.2in

The Segre-Verlinde equality tells us that expressions \eqref{segreexp1} and \eqref{segreexp3} are zero, as all terms of degree equal to $r^2(g-1) \mod 2$ in $\ell +1$ vanish in the Verlinde polynomial \eqref{verexp}. Expression \eqref{segreexp2} can in turn be matched to the much simpler term of degree $t+g -2$ in \eqref{verexp}.

We note that non-trivial identities are also obtained in {\it low degree} from the Segre-Verlinde comparison. In particular it is interesting to examine the known vanishing (cf. \cite{ek}) in degrees higher than $(r^2 -r) (g-1)$ of the Pontryagin ring of $\mathsf N (r, d)$. The latter is the cohomology subring generated by the $\mathsf a$ classes,  with our normalization  \eqref{normalize}. The vanishing implies that all terms of degree lower than 
$(r-1) (g-1)$ in $\ell +1$ are zero in the Verlinde polynomial \eqref{verlindeN}. We note below that when $r = 2,$ the Segre-Verlinde correspondence {\it{recovers}} in fact the Pontryagin vanishing. 

\vskip.2in

\subsection{The Segre-Verlinde identity in rank 2} 
We now turn to the calculation of the full Segre polynomial in the rank-two, odd-degree case. We write $\mathsf N_g (2, d)$ to emphasize the genus of the underlying curve of the moduli space. The latter has dimension $3g-3.$ It is convenient to write the Segre polynomial in terms of the formal variable $$\mathsf y = \frac{\alpha_1 - \alpha_2}{2}, \, \, \text{satisfying} \, \, \mathsf y^2 = - \mathsf a_2.$$
As before $\alpha_1, \alpha_2$ denote the Chern roots of the rank 2 universal bundle $\mathcal W_p \to \mathsf N_g (2, d).$ Starting with the Chern character expression \eqref{segreMN}, we are led to the formula
\begin{equation}
s(\alpha_{\mathsf M}) = (\ell +1)^g \int_{\mathsf N_g (2,\, d)} (1 + \mathsf a_2)^{g-2} \exp \left ( (\ell +1) \left [ \mathsf f_2 \cdot \frac{1}{\mathsf y} 
\ln \frac{1+\mathsf y}{1-\mathsf y} + \sum_{j=1}^g \mathsf b_2^j \mathsf b_2^{j+g} \cdot \frac{1}{2\mathsf y^3} \left [ \ln \frac{1+\mathsf y}{1-\mathsf y} - 2 \mathsf y \right ] \right ] \right ).
\label{segrerank2}
\end{equation}
On the Verlinde side, we have
\begin{equation}
 \chi (\mathsf M_g (2,  d), \, \Theta_2^\ell \otimes \Theta_1) = (\ell +1)^g \int_{\mathsf N_g (2, \,d)} \exp \, (2(\ell +1) \, \mathsf f_2) \cdot  \left ( \frac{\mathsf y}{\sinh \mathsf y} \right )^{2 (g-1)}.
\label{verlinderank2}
\end{equation}
The correspondence asserts that the two series \eqref{segrerank2} and \eqref{verlinderank2} evaluate identically on the fundamental class of the moduli space $\mathsf N_g (2, d).$ 
 
 To understand this equality, we set $\mathsf b =  \sum_{j=1}^g \mathsf b_2^j \mathsf b_2^{j+g} ,$ 
and note first the well-known identity, explained in \cite{thaddeus},
\begin{equation} \int_{\mathsf N_g (2, \,d)} \mathsf b^p \cdot \mathsf P(\mathsf a_2, \mathsf f_2)  =g \cdot  \int_{\mathsf N_{g-1}  (2, \,d)} \mathsf b^{p-1} \cdot \mathsf P (\mathsf a_2, \mathsf f_2),
\end{equation}
where $\mathsf P$ is any polynomial in the universal classes $\mathsf a_2, \mathsf f_2$ and $p < g.$
Assuming this identity, a careful analysis shows that the equality of the two integrals  \eqref{segrerank2} and \eqref{verlinderank2} for all genera gives the known formula (cf. \cite{thaddeus}), 
\begin{equation}
\int_{\mathsf N_g (2, \, d)} \mathsf f_2^m \mathsf a_2^n = \frac{m}{2} \int_{\mathsf N_{g-1} (2,\, d)} \mathsf f_2^{m-1} \mathsf a_2^{n-1}.
\label{rank2symmetry}
\end{equation}
(Note that relative to the classic Newstead notation \cite{newstead}, we have $\alpha = 2 \mathsf f_2, \, \beta = - 4 \mathsf a_2.$) 

Furthermore, we see that the Segre series \eqref{segrerank2} has no term of degree 0 in $\ell +1.$ For every odd genus, the correspondence therefore forces the vanishing of the degree $0$ term in the Verlinde series, giving
$$\int_{\mathsf N_g (2, d)} \mathsf a_2^{\frac{3g-3}{2}} =0.$$
Through the symmetry \eqref{rank2symmetry}, this establishes the general Pontryagin vanishing $$\int_{\mathsf N_g (2, \, d)} \mathsf f_2^m \mathsf a_2^n = 0, \, \, \, \text{whenever} \, \, \, m < n.$$
We see therefore that the Segre-Verlinde correspondence does not force the explicit evaluation of each top monomial $ \mathsf f_2^m \mathsf a_2^n,$ but comes down instead to the identity \eqref{rank2symmetry}, alongside the Pontryagin vanishing. The additional datum needed to pin down all top intersections is the evaluation of the volumes $\int_{\mathsf N_g (2, d)} \mathsf f_2^{3g-3}.$ Overall, we note that the emerging constraints bear similarities to the Virasoro constraints on $\mathsf N_g (2,d)$ described in \cite{blm}. 

Finally, we remark that in arbitrary rank, all $\mathsf f$ classes appear in the Segre polynomial in a factor of the form $$\exp{\left [ (\ell+1) (\mathsf f_2 \cdot P_2 (\mathsf a) +\,  \cdots \, + \mathsf f_r \cdot P_r (\mathsf a)) \right]},$$ This allows one to express intersections involving higher $\mathsf f$ classes in terms of the basic Verlinde intersections which involve only the classes $\mathsf f_2, \, \mathsf a_2, \ldots \mathsf a_r.$

\vskip.4in

\section{An intersection-theoretic triangle}

\vskip.2in

A relation between Segre and Verlinde invariants on moduli spaces of sheaves was first proposed in \cite{johnson} for the Hilbert scheme of points on a smooth projective surface. The correspondence of these seemingly different invariants was further pursued in \cite{mop1}, \cite{mop2}, \cite{gk}, \cite{gottsche}, \cite{oberdieck} over the Hilbert scheme of points as well as over higher-rank moduli spaces of stable sheaves on surfaces. In this context, a general proof of the Segre-Verlinde correspondence on the Hilbert scheme was given in \cite{gm}. A beautiful identity of Segre and Verlinde numbers was further uncovered in the case of Quot schemes of rank zero quotients over curves and surfaces in \cite{ajlop}. The correspondence was also investigated for Quot schemes of rank zero quotients on surfaces and fourfolds in \cite{bojko}, \cite{bh}. Theorem \ref{coprime} and its variant \ref{coprime}$^*$ place the 
classic moduli space of stable bundles over a curve into the emerging general framework of the Segre-Verlinde duality.

The current paper is not shaped however around a duality. At its core is the equality of three distinct invariants,
$$\chi (\M(r, d), \, \Theta_r^\ell \otimes {\det}^\star \Theta_1) =  \int_{[\Q (\mathbb C^{(\ell + h) r_0}, r, d)]^{vir}} a_r^{(\ell + h) d_0 - \ell r_0 (g-1)} =  \int_{\M (r, d)} s (\alpha_{\M}),$$
where $\alpha_{\M}$ corresponds to a $K$-class $\alpha$ on $C$ with $$\text{rank} \, \alpha = (\ell + h) r_0, \, \, \deg \alpha = -(\ell + h) d_0 + \ell r_0 (g-1).$$
The intersection number on the Quot scheme, important for the level-rank symmetry of the Verlinde numbers, plays a key role in the Segre-Verlinde identity formulated in Theorems \ref{coprime}-\ref{coprime}$^*$. In the surface case, the calculations and conjecture of \cite{johnson} also have as starting point the strange duality conjecture. Nevertheless, in a surface setting, a compelling formulation of the Quot scheme intersection in the triangle of invariants has remained elusive, and work surrounding it is in progress. It is natural as well to formulate the Segre-Verlinde correspondence for moduli spaces of Higgs bundles, a topic which will be addressed in future work.

\vskip.3in

\section*{Acknowledgements}

I thank Dragos Oprea and Rahul Pandharipande for many discussions surrounding the geometry of Quot schemes over the years. I am very grateful to the Mathematics Section of ICTP for making possible a wonderful stay in Trieste in the fall of 2021, when related ideas were explored. I thank Lothar G\"{o}ttsche for interesting conversations on Segre and Verlinde invariants.

\vskip.3in

\section*{Ethics declarations}
\subsection*{Competing interests} The author has no relevant financial or non-financial interests to disclose. 
\subsection*{Funding} This work was supported by the NSF through grant DMS 1902310. 

\vskip.3in

\section*{Data availability statement} Data sharing is not applicable to this article as no datasets were generated or analyzed during the current study.


\vskip.4in

\begin{thebibliography}{1} 

\bibitem [AB] {ab}
M. Atiyah, R. Bott, {\it The Yang-Mills equations over Riemann surfaces}, Phil. Trans. R. Soc. Lond. A {\bf 308} (1982), 523-615.

\bibitem [AJLOP] {ajlop}
N. Arbesfeld, D. Johnson, W. Lim, D. Oprea, R. Pandharipande, {\it The virtual K-theory of Quot schemes of surfaces,} J. Geom. Phys. {\bf 164} (2021). 

\bibitem [Be] {bertram}
A. Bertram, {\it Towards a Schubert calculus for maps from a Riemann surface to a Grassmannian,} Internat. J. Math. {\bf 5} (1994), no 6, 811-825.

\bibitem [Bo] {bojko}
A. Bojko, {\it Wall-crossing for punctual Quot schemes}, arXiv:2111.11102.


\bibitem [BH] {bh}
A. Bojko, J. Huang, {\it Equivariant Segre and Verlinde invariants for Quot schemes,} arXiv:2303.14266v2.

\bibitem [BLM] {blm} 
A. Bojko, W. Lim, M. Moreira, {\it  Virasoro constraints on moduli of sheaves and vertex algebras}, Invent. Math. {\bf 236} (2024), 387 - 476. 

\bibitem [Br] {bradlow}
S. Bradlow, {\it Special metrics and stability for holomorphic bundles with global sections}, J. Diff. Geom. {\bf 33} (1991), 169-214.

\bibitem [BrD]{bd}
S. Bradlow, G. Daskalopoulos, {\it Moduli of stable pairs for holomorphic bundles over Riemann surfaces,} Int. J. Math. {\bf 2} (1991), 477-513.

\bibitem [BDW] {bdw}
A. Bertram, G. Daskalopoulos, R. Wentworth, {\it Gromov Invariants for holomorphic maps from Riemann surfaces to Grassmannians,} J. Amer. Math. Soc. {\bf 9} (1996), no 2, 529-571.

\bibitem[DN]{dn}
J.-M. Dr\'{e}zet, M. Narasimhan, {\it Groupe de Picard des vari\'{e}t\'{e}s de modules de fibr\'{e}s semistables sur les courbes alg\'{e}briques}, Invent. Math. {\bf 97} (1989), 53-94.

\bibitem [EK] {ek}
R. Earl, F. Kirwan, {\it The Pontryagin rings of moduli spaces of arbitrary rank holomorphic bundles over a Riemann surface,} J. London Math. Soc. (2) {\bf 60} (1999), 835-846. 

\bibitem [G] {gottsche}
L. G\"{o}ttsche, {\it Blowup formulas for Segre and Verlinde numbers of surfaces and higher rank Donaldson invariants,} Proc. Sympos. Pure Math. {\bf 109} (2024), 127 - 152. 

\bibitem [GK] {gk}
L. G\"{o}ttsche, M. Kool, {\it Virtual Segre and Verlinde numbers of projective surfaces}, J. Lond. Math. Soc. (2) {\bf 106} (2022), 2562 - 2608. 

\bibitem [GM] {gm}
L. G\"{o}ttsche, A. Mellit, {\it Refined Verlinde and Segre formula for Hilbert schemes}, arXiv:2210.01059. 

\bibitem [I]{i}
K. Intriligator, {\it Fusion Residues,} Modern Phys. Lett. A {\bf 6} (1991), 3543-3556.

\bibitem [J] {johnson}
D. Johnson, {\it Universal series for Hilbert schemes and strange duality}, IMRN {\bf 10} (2020) 3130–3152.

\bibitem[L] {lin}
Y. Lin, {\it The moduli space of stable pairs,} Pacific J. Math. {\bf 294} (2018), no. 1, 123-158.

\bibitem [M] {alina}
A. Marian, {\it On the intersection theory of Quot schemes and moduli of bundles with sections,} J. reine angew. Math. {\bf 610} (2007), 13-27.

\bibitem[MO1]{mo1}
A. Marian, D. Oprea, {\it Virtual intersections on the Quot scheme and Vafa-Intriligator formulas,} Duke Math. J. {\bf 136} (2007), 81-113.


\bibitem [MO2]{mo2}
A. Marian, D. Oprea, {\it The level-rank duality for non-abelian theta functions,}  Invent. Math. {\bf 168} (2007), 225 - 247.

\bibitem [MO3]{mo3}
A. Marian, D. Oprea, {\it On the intersection theory of the moduli space of rank two bundles,} Topology {\bf 45} (2006), 531-541.

\bibitem [MO4]{mo4}
A. Marian, D. Oprea, {\it Counts of maps to Grassmannians and intersections on the moduli space of bundles,} J. Diff. Geom. {\bf 76} (2007), 155-175. 

\bibitem [MOP1] {mop1}
A. Marian, D. Oprea, R. Pandharipande, {\it The combinatorics of Lehn’s conjecture}, J. Math. Soc. Japan {\bf 71} (2019), 299–308.

\bibitem [MOP2] {mop2}

A. Marian, D. Oprea, R. Pandharipande, {\it Higher rank Segre integrals over the Hilbert scheme of points}, J. Eur. Math. Soc. {\bf 24} (2022), 2979–3015.

\bibitem [N] {newstead} 

P. Newstead, {\it Characteristic classes of stable bundles of rank $2$ over an algebraic curve,} Trans. Amer. Math. Soc. {\bf 169} (1972), 337 - 345.  

\bibitem[O]{oberdieck}
G. Oberdieck, {\it Universality of descendent integrals over moduli spaces of stable sheaves on $K3$ surfaces,} SIGMA {\bf 18} (2022). 

\bibitem[R] {r}

K. Rietsch, {\it Quantum cohomology rings of Grassmannians and total positivity}, Duke Math. J. {\bf 110} (2001), 523 - 553.

\bibitem [ST] {st}
B. Siebert, G. Tian, {\it On quantum cohomology rings of Fano manifolds and a formula of Vafa and Intriligator,} Asian J. Math. {\bf 1} (1997). 679-695. 

\bibitem [T] {thaddeus} 
M. Thaddeus, {\it An introduction to the topology of the moduli space of stable bundles on a Riemann surface,} Lecture Notes in Pure and Appl. Math. {\bf 184}, 71 - 99. 

\bibitem [W] {witten}
E. Witten, {\it The Verlinde algebra and the cohomology of the Grassmannian,} Conf. Proc. Lect. Notes Geom. Topol. vol. 4, p. 357-422. Internat. Press, Cambridge, MA (1995). 


\end{thebibliography}
\end{document}